\documentclass[11pt]{amsart}          	
\usepackage{amssymb,latexsym,amsmath,amscd,stmaryrd}
\usepackage{algorithm}
\usepackage{algorithmic}

\newcommand{\F}{{\mathbb{F}}}
\newcommand{\C}{{\mathbb{C}}}

\newcommand{\ds}{\displaystyle}

\newcommand{\m}{\mathfrak{m}}

\newcommand{\CO}{{\mathcal O}}

\newcommand{\CR}{{\mathcal R}}
\newcommand{\CK}{{\mathcal K}}

\newcommand{\ARA}{\begin{array}}
\newcommand{\ARE}{\end{array}}
\newcommand{\SMA}{\begin{smallmatrix}}
\newcommand{\SME}{\end{smallmatrix}}
\newcommand{\BMA}{\begin{matrix}}
\newcommand{\BME}{\end{matrix}}

\newtheorem{lem}{{Lemma}}[section]

\newtheorem{prop}[lem]{{Proposition}}

\newtheorem{defn}[lem]{{Definition}} 
\newtheorem{hypo}[lem]{{Hypothesis}} 
\newtheorem{rem}[lem]{{Remark}} 

\theoremstyle{definition}
\newtheorem{exmp}[lem]{{Example}} 
\newcommand{\lra}{\longrightarrow}

\newcommand{\ra}{\rightarrow}

\DeclareMathOperator{\NF}{NF}

\DeclareMathOperator{\Ann}{Ann}
\DeclareMathOperator{\spec}{Spec}
\DeclareMathOperator{\Sing}{Sing}
\DeclareMathOperator{\iso}{\widetilde{\rightarrow}}
\DeclareMathOperator{\embdim}{embdim}
\DeclareMathOperator{\minor}{minor}
\begin{document}
\title{Milnor algebras could be isomorphic to modular algebras}
\author{Bernd Martin, Hendrik S\"u\ss}




\newpage
\setcounter{page}{1}
\maketitle
\begin{abstract}
  We find and describe unexpected isomorphisms between two very 
  different objects associated to hypersurface singularities.
  One object is the Milnor algebra of a function, while the other object 
  is the local ring of the flatness stratum of the 
  singular locus in a miniversal deformation, an invariant of the
  contact class of a defining function. 
  Such isomorphisms exist for unimodal hypersurface singularities. 
  However, for the moment it is not well understood, which principle causes these 
  isomorphisms and how far this observation generalises. 
  Here, we also provide  
  an algorithmic approach for checking algebra isomorphy.  
\end{abstract}

\section{Introduction}
Let $X_0\subseteq \C^n$ be a germ of an isolated hypersurface 
singularity  defined by an analytic function 
$f(x)=0$, $ f\in \C\{x\}$. An important topological 
invariant of the germ is the Milnor number, which can be computed as 
the $\C$-dimension
of the so-called Milnor algebra 
$Q(f)=\C\{x\}/(\partial f/\partial x)$, \cite{1}.
The Milnor algebra carries a canonical structure of a 
$\C[T]$-algebra defined by the multiplication with $f$. 
A special version of the Mather-Yau theorem states that the $\CR$-class 
(right-equivalence class) of the function $f(x)$ with isolated critical point 
is fully determined by the isomorphism class of 
$Q(f)$ as $\C[T]$-algebra \cite{scherk}.

By computational experiments, we have found another occurrence of the 
Milnor algebra -- this time connected with the $\CK$-class 
(the contact equivalence class) of $f(x)$, i.e. with the   
isomorphism-class 
of the germ $X_0$. Our observation concerns 
unimodal functions that are not quasihomogeneous.
Here we consider a miniversal deformation $F:X\ra S$ of the 
singularity $X_0$. 
It has a smooth base space of dimension $\tau$, with $\tau$ being the Tjurina number, 
i.e. the $\C$-dimension of the Tjurina algebra $T(f):=Q(f)/fQ(f)$.
We consider the relative singular locus $\Sing(X/S)$ of $X$ over $S$ 
and its flatness stratum $\F:=\F_S(\Sing(X/S))\subset S$, 
which depends only on $X_0$, up to isomorphism. 
The flatness stratum is computable for
sufficiently simple functions using a special algorithm \cite{flatalg}.
Surprisingly, the local ring of the flatness stratum  of a unimodal
singularity is isomorphic in all computed cases, either 
to the Milnor algebra
of the defining function (in case $\dim(\F)=0$), 
or to the Milnor algebra of a 'nearby'
function with non-isolated critical point, otherwise.  

The notion of a modular stratum was developed by Palamodov, 
\cite{P1}, in order to find a moduli space for singularities. 
It coincides with the 
flatness stratum $\F=\F_S(\Sing(X/S))$ \cite{flat}, which has been described 
for unimodal functions in \cite{lum}. 
Only for some singularities from the T-series the modular stratum has 
expected dimension 1 with smooth curves and embedded fat points as 
primary components. The combinatorial pattern of its occurrence 
was found and the phenomenon of a splitting singular locus along a 
$\tau$-constant stratum was discovered. Here we extend our 
observation that the modular stratum 
is the spectrum of the Milnor algebra of an associated non-isolated 
limiting singularity.

The modular stratum is a fat point of multiplicity $\mu$
isomorphic to $\spec Q(f)$ in all other (computed) cases 
of $T$-series singularities.
The same holds for all 14 
exceptional and non-quasihomogeneous unimodal singularities.
In the case of a quasihomogenous exceptional singularity, the modular stratum 
is a smooth germ, hence corresponding to a trivial Milnor algebra. 
 
For completeness, we will first recall the basic results on modular strata 
and prove that they are algebraic. 
Second, we collect and complete results
on the modular strata of unimodal functions, which are already
found in \cite{lum}. 
Subsequently, some of the non-trivial unexpected isomorphisms are presented.
A further example of higher modality is discussed in section 4. 
Hypotheses towards a possible 
generalisation of these experimental
results are formulated. 
Finally an algorithmic approach is outlined, how to perform a check of algebra
isomorphy using a computer algebra system, 
knowing that there is no practical algorithm in general.
All computations were executed in the computer algebra system 
{\sc Singular} \cite{9}.   

\section{Characterisations of a modular germ}
The definition of modularity was introduced by Palamodov, cf.\  
for instance \cite{P1} and \cite{P2}, and was
simultaneously discussed by Laudal for the case of formal 
power series under the name {\em prorepresentable} substratum'.
This notion can be considered for any isolated singularity with 
respect to several deformation functors or to deformations
of other objects, cf.~\cite{nato}. For simplicity we restrict ourselves to the case of
a germ of an isolated complex hypersurface singularity 
$X_0=f^{-1}(0)\subseteq \C^n$, or an isolated complete intersection 
singularity (ICIS).

A {\em deformation} of $X_0$  is a flat 
morphism of germs $\xi:X\lra S$ with its  special fibre isomorphic to $X_0$. 
It is called {\em versal}, if any other deformation of $X_0$ can be 
induced via a morphism of the base spaces up to isomorphism. 
It is called {\em miniversal}, if the dimension of the base
space is minimal. Miniversal deformations exist for isolated 
singularities and are unique up to a non-canonical isomorphism.
In case of a hypersurface, a miniversal deformation has a smooth base 
space, i.e. the deformations are unobstructed. 
It can be represented as an {\em embedded} 
deformation
$\xi:X\subset \C^n\times S \ra S$,
$S=\C^\tau$, $X=F^{-1}(0)$, $F(x,s)=f(x)+
\sum_{i=1}^\tau s_\alpha m_\alpha$, 
where $\{m_1,\ldots , m_\tau\}\subset\C\{x\}$ 
induces a $\C$-basis of the Tjurina algebra $T(f)$. 

Obviously, a miniversal deformation has not the properties of a 
moduli space, because there are always isomorphic fibres or 
even locally trivial subfamilies.
Hence, the inducing morphism of another deformation is not unique. 
One can, however, look for subfamilies of a miniversal deformation having
this universal property. 
 
\begin{defn}
Let $\xi:X\ra S$ be a miniversal deformation of a complex germ $X_0$. A subgerm 
$M\subseteq S$ of the base space germ is called {\em modular} if the 
following universal property holds: If $\varphi:T \ra M$ 
and $\psi:T\ra S$ are morphisms such that the induced deformations 
$\varphi^*(\xi_{|M})$ and $\psi^*(\xi)$ over $T$ are isomorphic, 
then $\varphi=\psi$.
\end{defn}
The union of two modular subgerms inside a miniversal family is again
modular. Hence, a unique maximal modular subgerm exists. 
It is called {\em modular stratum} of the singularity. Note that  
any two modular strata of a singularity are isomorphic by definition.
\begin{exmp} \label{ex1}
If $X_0$ is an ICIS with a good
$\C^*$-action, i.e. defined by quasihomogeneous  polynomials, 
then its modular stratum coincides with the reduced $\tau$-constant stratum and is smooth, 
cf.~\cite{alex}.
\end{exmp}
Palamodov's definition of modularity
is difficult to handle.  It made it challenging
to find non-trivial explicit examples. 
Even the knowledge of the basic characterisations 
of modularity in terms of cotangent cohomology, which were  already discussed 
by Palamodov and Laudal, did 
not lead to more examples, cf.~\cite[Thm. 6.2]{P1}.
   
\begin{prop} 
Given a miniversal deformation $\xi:X\ra S$
of an isolated singularity $X_0$,
the following conditions are equivalent for a subgerm of the base space 
$M\subseteq S$:
\begin{enumerate}
\item $M$ is {\em modular}.   
\item $M$ is {\em infinitesimally modular}, i.e. injectivity of the 
      relative Kodaira-Spencer map 
      $T^0(S,\CO_M)$ $ \lra$ $T^1(X/S,\CO_S)_{|M}$ holds.
\item $M$ has the {\em lifting property} of vector fields of 
      the special fibre:
      $T^0(X/S,\CO_S)_{|M}$ $\lra$ $T^0(X_0,\C)_{|M}$ is surjective.       
\end{enumerate} 
\end{prop}  
Note that $T^0$ corresponds to the module of associated vector fields,
while $T^1$ describes all infinitesimal deformations. It is 
represented in the hypersurface case by the (relative) Tjurina algebra 
$T^1(X/S)=T(F)=\C\{x,s\}/(F,\partial_x F)$. 

All geometric objects belong to the category of analytic germs. 
But an isolated singularity is algebraic, i.e. its defining
equations can be chosen as polynomials.
Ad hoc it is not clear whether the modular stratum is algebraic, too,
and to our knowledge it has not been demonstrated.
Here, we add a proof in case of ICIS for completeness.

\begin{lem}
Let $X_0$ be a germ of an isolated complete intersection singularity.
Then its modular stratum $M(X_0)\subset\C^\tau$ is an algebraic subgerm.
\end{lem}
The proof uses the characterisation of modularity as flatness stratum
of the Tjurina-module. A more general result holds under 
weaker assumptions than ICIS, too,  cf.~\cite[Prop. 2.1]{flat}.
\begin{prop}
\mbox{} \\
 Let $X_0\subseteq \C^n$ be an isolated complete intersection 
 singularity defined by $p$ equations $f\in \C\{x\}^p$ 
 with miniversal deformation
 $\xi:X\ra S$. Then the modular space coincides with the flatness stratum
 of the relative Tjurina module 
 $T^1(X/S)=\CO^p_{X}/(\partial F/\partial x)\CO^p_{X}$ as 
$\CO_S$-module, $F$ being the equations of the deformation.
\end{prop}

\begin{proof}[Proof of the lemma]
We may choose the defining equations $f=(f_1, \ldots , f_p)$ of the germ $X_0$ 
as polynomials 
by finite determination of isolated singularities. 
The affine variety defined by
these polynomials $V(f)\subset \C^n$ has in general other singularities than the zero point.
But, we can choose an embedding such that 
$\Sing(V(f))$ is concentrated at zero. This holds iff global and local 
Tjurina number are equal 
$$\dim_\C(\C[x]^p/(f\C[x]^p,\partial f/\partial x))= \dim_\C(\C\{x\}^p/(f\C\{x\}^p,\partial f/\partial x))=\tau.$$
Consider the $\C[s,x]$-module 
$B:=\C[s,x]^p/(F\C[s,x]^p,\partial F/\partial x)$. The module $B$ is finite as $\C[s]$-module.
Its flatness stratum over $S$ at zero $\F_{S,0}(B)\subset S$, $S:=\C^\tau=\spec(\C[s])$,  
is well defined by the fitting ideal of $B$ as $\C[s]$-module.
The $\C\{s,x\}$-module $T^1(X/S,\CO_S)$ is finite as $\C\{s\}$-module.
Consider the modules $B_0:=B/sB$ and $T^1(X_0)=T^1(X/S,\CO_S)_{|s=0}$, then  the
localisation at $x=0$ of $B_0$ and $T^1(X_0)$ have identical module-structures which are 
both already given as $\C[x]/(x)^k$-modules: $B_{0\ (x)}=T^1(X_0)$, hence the germ at zero
$\F(B)_{(s,x)}$ coincides with the flattening stratum of $T^1(X/S,\CO_S)$.
\end{proof}

We add some remarks concerning the flatness criterion:
\begin{itemize}
\item 
The support of $T^1(X/S,\CO_S)$ is exactly the relative singular locus of the
mapping germ $F:\C^n\times S \lra \C^p\times S$ over $S$. In case of a
hypersurface, i.e. $p=1$, $T^1(X/S,\CO_S)$ coincides with the $\CO_S$-algebra of the 
relative singular locus, that is the relative Tjurina-algebra $T(F)=\CO_{\Sing(X/S)}$. 
\item The support of the 
flatness-stratum $\F_{\CO_S}(\Sing(X/S))$
is 
the locus where the finitely generated $\CO_S$-module $T^1(X/S,\CO_S)$ has constant rank $\tau$.
The fact that $M_{red}$ is the $\tau$-constant stratum was already explained by Palamodov, cf.~\cite[Thm. 7.1]{P1}.
Besides the case of a quasihomogeneous singularity, see example \ref{ex1}, a non-reduced structure 
of $M$ is expected generally.
\item It follows from a non-trivial result, cf.~\cite{le}, that the finite map 
$\Sing(X/S) \rightarrow S$ is unramified  over the $\mu$-constant stratum. 
But the analogous statement for the $\tau$-constant stratum 
does not hold, see below. This phenomenon we have called {\em splitting singular locus}
over the $\tau$-constant stratum.
\end{itemize}
\section{Computing the modular germs of unimodal singularities}
Applying the algorithm for computing the flatness stratum, cf.~\cite{flatalg}, we can compute the modular 
stratum of 
not too complicated 
singularities. More precisely, the output
of the algorithm is the $k$-jet of the germ of the flatness
stratum for some positive integer $k$. If the modular stratum is a fat point 
we are done with some big number $k$. We 
can not prove or even expect to end up with an algebraic representation generally. 
But, it does occur, as in all presented examples below.

The classification of singularities starts with the simple
singularities, the ADE-singula\-rities. These are all quasihomogeneous,
their modular strata are all trivial, i.e. simple points.
Following the classification of functions by Arnol'd \cite{agv}
the next more complicated singularities are the unimodal ones. They are 
characterised by the fact, that in a neighbourhood of the function
only $\CR$-orbit families occur which depend on at most one parameter.  
 
Recall the classification of unimodal functions: 
We have the $T$-series singularities and 14 so called 
exceptional unimodal singularities. 
We may restrict their representation 
to three variables up to stable equivalence.
Any type is representing 
a one-parameter $\mu$-constant family of $\CR$-equivalence 
classes.
The exceptional ones, cf.~appendix~\ref{sec:appendix-table}, are all semi-quasihomogeneous of a special type:
The $\mu$-constant family can be written as
$$f_\lambda=f_0(x)+\lambda h_f(x),\ \ \lambda\in\C,
$$ 
where $f_0$ is quasihomogeneous and 
$h_f(x):=\det(\frac{\partial^2f_0}{\partial x_i \partial x_j})$ is the Hesse form of 
$f_0$. Such a family splits into exactly  two $\CK$-classes,
one quasihomogeneous ($\lambda=0$) with trivial modular stratum and one semi-quasihomogeneous  
($\lambda \not= 0$, we call it of {\em Hesse-type}). Here $\tau(f_1)=\mu(f_1)-1$ holds, the 
modular strata are fat points of multiplicity $\mu$, cf.~\cite[Prop. 5.1]{lum}.  
     
The singularities of the $T$-series are defined by the equations
$$T_{p,q,r}: \ x^p+y^q+z^r+\lambda xyz, 
\ \ \ds{\frac{1}{p}+\frac{1}{q}+\frac{1}{r}\leq 1}.$$
 
In exactly three cases we have 
$ \ds{\frac{1}{p}+\frac{1}{q}+\frac{1}{r}=1}$.
These singularities are quasihomogeneous. 
They are called the {\em parabolic} singularities 
$P_8$, $X_9$ and $J_{10}$
in Arnold's notation or {\em elliptic} hypersurface singularities 
$\tilde{E}_6$, $\tilde{E}_7$ and $\tilde{E}_8$, 
in Saito's paper \cite{sa}:
 
\begin{itemize}
\item[] \ \ $\tilde{E}_6=P_8= T_{3,3,3}: \ x^3+y^3+z^3+\lambda xyz, \;\lambda^{3}\neq-3^{3},\;\tau=\mu=8$; 
\item[] \ \ $\tilde{E}_7=X_9=T_{4,4,2}: \ x^4+y^4+z^2+\lambda xyz, \;\lambda^{4}\neq2^{6},\;\tau=\mu=9$; 
\item[] \ \ $\tilde{E}_8=J_{10}=T_{6,3,2}: \  x^6+y^3+z^2+\lambda xyz, \;\lambda^{6}\neq2^{4}3^{3},\;\tau=\mu=10$.
\end{itemize}
It is well known that their $\CK$-classes 
(here equal to the $\CR$-classes)
are not determined by $\lambda$. The
$\CK$-equivalence induces a discrete equivalence relation on the 
$\lambda$-line with the indicated gaps. 
Its quotient is an affine line parametrised by
the classical $j$-invariant of elliptic curves. 
Precise formul\ae\   for $j$ can be found in Saito's paper. Their modular strata are germs of a line.

All other $T$-series singularities are called  
{\em hyperbolic}. Their $\CK$-class is independent of $\lambda$, $\lambda\not=0$, the 
Newton boundary has three maximal faces and the singularity is
neither quasihomogeneous nor semi-quasihomogeneous. 
We have $\tau(T_{p,q,r})=\mu(T_{p,q,r})-1=p+q+r-2$.
The modular strata of the hyperbolic singularities are more 
complicated. Some of them are 1-dimensional, others are just fat points.
First computations are discussed  in \cite{lum}. 

Obviously, a hyperbolic singularity of type $T_{p,q,r}$ is adjacent to another T-series
singularity, iff all three parameters ($p,q,r$) are greater or equal to the parameters
of the second. Hence, any hyperbolic singularity is adjacent to at least one 
parabolic singularity. 
Inspecting the list we find exactly six types of
hyperbolic singularities which have the same Tjurina number as an adjacent parabolic
singularity. They are candidates for modular strata of dimension 1, because in their miniversal deformation 
the associated parabolic singularity occurs as one dimensional $\tau$-constant family including the special fibre.
Indeed, these six singularities are leading members of six sub-series of the hyperbolic singularities, which have 
a one-dimensional modular strata, cf.~\cite[Prop.4.2, 4.3]{lum}. These exceptional sub-series are characterised by the 
fact, that two of the three indices coincide with that of a parabolic one.

\begin{prop}
Any singularity of one of the following six exceptional sub-series contains a splitting line in their modular stratum:
 $$(T_{k,3,3})_k,\ k\geq l=4,\ \ \ \  (T_{k,4,2})_k,\ k\geq l=5,\ \ \ \ (T_{k,4,4})_k, \ k\geq l=3,$$
  $$(T_{k,3,2})_k,\ k\geq l=7, \ \ \ \ (T_{k,6,2})_k,\ k\geq l=4, \ \ \ \ (T_{k,6,3})_k, \ k\geq l=3.$$
\end{prop}
The families over the $\tau$-constant lines with index $k$ are given by 
  $$ f_t:=x^{l-1}(x+t)^{k-l+1}+y^q+z^r+xyz.$$  
  The fibre singularities over $t\not=0$ consist of one singularity of 
the associated parabolic type and of one singularity of type $A_{k-l-1}$ outside the zero section, if $k>l-1$.
These lines, called {\em splitting lines}, are components of the modular stratum, which
mostly contains another embedded component (fat point) at zero.
Exactly three of the singularities of the sub-series have two line components and two of them
have three line components, cf.~\cite[Cor. 4.4]{lum}. 
We call these five types the {\em symmetric exceptions}, because they come
from an additional symmetry in the equation or from a crossover of two exceptional sub-series.
All other computed examples of modular strata of $T$-series 
singularities, not belonging to the above
six exceptional sub-series, are fat points. 
\begin{exmp}
$f_t:=x^4+y^3+z^3+xyz+tx^3$ is a $\tau$-constant  deformation of $T_{4,3,3}$ with
generic fibre type $P_8$.
The modular deformation $f_t$ fits into the $\lambda$-line 
of $P_8(\lambda)$  at infinity: 
$f_t\sim_\CK P_8(t^{-1/3})$ for $t\not= 0$, i.e. we get a threefold covering 
of the $t$-line, $t\not=0$, by the $\lambda$-line, $\lambda\not=0$. 
We may think of a compactification of the modular $\lambda$-line of $P_8$ at 
infinity with a point corresponding to the $T_{4,3,3}$-singularity.
The same holds for $T_{4,4,3}$ and $T_{5,4,2}$ with respect to $X_9$.
This causes two different compactifications of the same modular family over the
punctured disc $X_9(1/\lambda)$, $\lambda>N$, at the special point zero 
to a modular family over the disc.
\end{exmp}

\section{New explicit results on unimodal modular strata}
A careful inspection of our
computations leads us to generic formul\ae\   for the equations of the modular strata 
of the $T$-series singularities in terms
of the index-parameters $(p,q,r)$. 
We confirm these equations by checking all modular strata
up to $\tau<50$, giving a strong indication of their correctness. Moreover, it gives
an explanation for the occurrence of line components in the modular strata of 
exceptional sub-series singularities.
In addition we obtain
the isomorphy of the modular algebra of a $T$-series singularity, not belonging to
one of the six sub-series, to its Milnor algebra,
the isomorphy of the modular algebra of sub-series singularities to the Milnor algebra of
certain non-isolated functions, and
the isomorphy of all modular strata from the same sub-series up to
the five symmetric exceptions.

\begin{exmp}
\label{exmp:tseries}
  Let $X_0$ be the hyperbolic singularity defined by $f=x^p+y^q+z^r +xyz$.  Then
  $$F=f + \mathbf{t_1} x^{p-1} + \ldots + \mathbf{t_{p-1}}x + \mathbf{t_p} + \mathbf{u_1} y^{q-1} + \ldots + \mathbf{u_{q-1}}y + \mathbf{v_1} z^{r-1} + \ldots + \mathbf{v_{r-1}}z$$
  defines a miniversal deformation $X \rightarrow S$ of $X_0$, with $\CO_S=\C\{t,u,v\}$.

  We consider the following the ideal $I(p,q,r) \subset \CO_S$: 
  \begin{eqnarray*}
    I(p,q,r) & = & (f_2,\ldots,f_{p}, g_2,\ldots,g_{q-1}, h_2,\ldots,h_{r-1},\\
    & &\; \mathbf{u_1v_1}-c_p \mathbf{t_1^{p-1}},\quad 
     \mathbf{t_1v_1}-d_q \mathbf{u_1^{q-1}},\quad 
     \mathbf{t_1u_1}-e_r \mathbf{v_1^{r-1}})
  \end{eqnarray*}
  where
     $$
     f_i := a^2\mathbf{t_{i}} - c_{i} \mathbf{t_1^i},\quad 
     g_i   := a^2\mathbf{u_{i}} - d_{i} \mathbf{u_1^i},\quad 
     h_i   := a^2\mathbf{v_{i}} - e_{i} \mathbf{v_1^i}
     $$
with coefficients
$$a:=a(p,q,r) := pqr(1-\frac{1}{p}-\frac{1}{q}-\frac{1}{r}),\ 
c_{i}:=\frac{\prod_{k=1}^{i}  a(p-k+1,q,r)}{i!a^{i-2}},$$
$$d_{i}:=\frac{\prod_{k=1}^{i}  a(p,q-k+1,r)}{i!a^{i-2}}, \
e_{i}:=\frac{\prod_{k=1}^{i}  a(p,q,r-k+1)}{i!a^{i-2}}.$$

The coefficient $c_p$ is zero iff  
$\frac{1}{k}+\frac{1}{q}+\frac{1}{r}=1$ for some $1 \leq k \leq p$, and similarly for $d_q$ and $e_r$.
So one of them vanishes exactly if $T_{p,q,r}$ belongs to one of the six exceptional sub-series.
More than one of the coefficients  is zero exactly for the five symmetric exceptions, see proposition \ref{symexc}.
Hence, the number of vanishing coefficients corresponds to 
the number of line components.
\end{exmp}
We have checked by computer.
\begin{prop}
\label{sec:prop-t-series-formulae}
The ideal $I(p,q,r)$ from example~\ref{exmp:tseries} defines 
the modular stratum of the singularity $T=T_{p,q,r}$, if $\tau(T) \leq 50$. 
\end{prop}
A linear diagonal transformation 
$$t_1 \mapsto \alpha t_1, u_1 \mapsto \beta u_1, v_1 \mapsto \gamma v_1.$$
induces isomorphy of the algebra defined by $I(p,q,r)$ and the Milnor algebra of $T_{p,q,r}$
in the case that none of the coefficients $c_p$, $d_q$ and $e_r$ vanishes. 
\begin{prop}
If $T_{p,q,r}$ does not belong to one of the six exceptional sub-series, then
$$Q(T_{p,q,r})\cong \CO/I(p,q,r).$$
In particular, for those singularities with $\tau<50$ the Milnor algebra is isomorphic to the
modular algebra of the singularity.
\end{prop}
We associate to an exceptional sub-series $(T_{k,q,r})_k$ the limit singularity $T_{\infty,q,r}$
defined by the equation $y^q+z^r+xyz$, which is non-isolated and belongs to the closure of the $\CK$-orbit
of any member of the sub-series.
\begin{prop}
If $T_\bullet=T_{k,q,r}$ belongs to one of the six exceptional sub-series, then
$$Q(T_{\infty,q,r})\cong \CO/I(k,q,r),$$
apart from the five symmetric exceptions.  
In particular, for those singularities $T_\bullet$ with $\tau<50$ the Milnor algebra  $Q(T_{\infty,q,r})$ is isomorphic to the
modular algebra of the singularity.
\end{prop}
Again, the isomorphy to the Milnor algebra is obtained 
by multiplying $t_1,u_1,v_1$ with some constants. It remains to check the five symmetric exceptions.
\begin{prop}\label{symexc}
 The local algebra of a
modular stratum is 
isomorphic to the Milnor algebra of a non-isolated singularity 
for the five symmetric exceptions of $T$-series singularities: 
\begin{eqnarray*}
Q(xyz)     & \mbox{for} \ \ T_{4,4,4}, & \mbox{and} \ \ T_{6,3,3},\\ 
Q(x^2+xyz) & \mbox{for} \ \ T_{6,4,2}  & \mbox{and} \ \ T_{6,6,2},\\ 
Q(x^3+xyz) & \mbox{for} \ \ T_{6,6,3}.
\end{eqnarray*}
\end{prop}
This follows from the equations in example~\ref{exmp:tseries} and proposition~\ref{sec:prop-t-series-formulae}.

Next we investigate the modular strata of all 14 exceptional 
semi-quasihomogeneous singularities and obtain a similar result.

\begin{prop}
All 14 exceptional semi-quasihomogeneous unimodal singularities fulfil:
 The local ring of their modular stratum is isomorphic to their Milnor algebra. 
\end{prop}

We list our results in a table at the end of the article. The isomorphisms are 
rather complicated. They are computed with the algorithm presented in appendix~\ref{sec:algorithm}.
We omit the formul\ae, but discuss only one case in detail.

\begin{exmp}[$W_{12}$, $Z_{11}$ and $S_{11}$]
  We consider $f = x^4+y^5+x^2y^3$  and choose 
  $(b_{11},\ldots, b_{1}):=(1,x,x^2,y,xy,x^2y,y^2,xy^2,xy^3,y^4)$ 
  as representatives of a $\C$-basis of the Tjurina algebra $T(f)$. 
  Now, $F=f + s_1b_1 + \ldots + s_{11}b_{11} \in \C\{x,y\} \otimes \CO_S$ defines 
  a miniversal deformation $X \rightarrow S$ of $X_0$.

 The ideal $I_M\subset \CO_S$ of the maximal modular subgerm 
 $M \subset S$, computed with {\sc Singular} is given 
 by the following completely interreduced generators:

  \begin{eqnarray*}
    \mathbf{s_{1}^4}&\mathbf{-}&\mathbf{\frac{30445}{7392}s_{1}^2s_{2}^2+\frac{4240139}{1897280}s_{1}^3s_{2}^2},\\
    \mathbf{s_{2}^3}&\mathbf{-}&\mathbf{\frac{2696}{48125}s_{1}^3s_{2}},\\
    s_{11}&+&\frac{11699}{144375}s_{1}^3s_2^2,\\
    s_{10}&-&\frac{3904}{48125}s_{1}^3s_2,\\
    s_9&+&\frac{52}{625}s_{1}^3-\frac{951}{7000}s_{1}s_2^2+\frac{592717}{8421875}s_{1}^2s_2^2-\frac{119567878949}{5187875000000}s_{1}^3s_2^2,\\
    s_8&+&\frac{1304}{5775}s_{1}^2s_2^2-\frac{1411481}{18528125}s_{1}^3s_2^2,\\
    s_7&-&\frac{618}{1925}s_{1}^2s_2+\frac{1024869}{37056250}s_{1}^3s_2,\\ 
    s_6&+&\frac{6}{25}s_{1}^2+\frac{3}{80}s_2^2-\frac{21}{3125}s_{1}^3+\frac{531}{20000}s_{1}s_2^2-\frac{31001023}{5390000000}s_{1}^2s_2^2,\\
    &+&\frac{25063327841}{207515000000000}s_{1}^3s_2^2,\\
  \end{eqnarray*}
  \begin{eqnarray*} 
    s_5&-&\frac{2}{25}s_{1}^3+\frac{9}{16}s_{1}s_2^2-\frac{114057}{539000}s_{1}^2s_2^2+\frac{6306416817}{83006000000}s_{1}^3s_2^2,\\
    s_4&-&\frac{6}{7}s_{1}s_2+\frac{1227}{67375}s_{1}^2s_2-\frac{16557777}{2593937500}s_{1}^3s_2,\\
    s_3&-&\frac{2}{5}s_{1}^2+\frac{9}{16}s_2^2-\frac{9}{625}s_{1}^3-\frac{621}{4000}s_{1}s_2^2+\frac{49325643}{1078000000}s_{1}^2s_2^2,\\
    &-&\frac{644553838881}{41503000000000}s_{1}^3s_2^2.    
  \end{eqnarray*}
$\CO_M$ is a zero-dimensional local algebra of 
 embedding dimension $2$. 
 A minimal embedding is defined by the two polynomials printed in bold.
The mapping   
\begin{eqnarray*}
\varphi:\CO_M &\rightarrow& \C\{x,y\}/(\textstyle \frac{\partial f}{\partial x}, \frac{\partial f}{\partial y})\\
\overline{s_{1}} &\mapsto& \frac{2668050}{2051993}\overline{y}-\frac{11759762521878525}{25638801731506361}\overline{y}^2\\
\overline{s_2} &\mapsto& \frac{2134440}{2051993}\sqrt{-\frac{1386}{6089}}\cdot\overline{x}
\end{eqnarray*}
defines an isomorphism between this local algebra and the Milnor algebra of $f$. 

We give the isomorphism for $Z_{11}$ and $S_{11}$:\\[1mm]
\begin{tabular}{ll}
name: &$S_{11}$\\
equation: &$f=y^2z+xz^2+x^4+x^3z$\\
deformation: & $F=f+s_{1} x^2z + s_2x^2y + s_3x^3 + s_4xz + s_5z +s_6xy+ s_7y + s_8x^2 + s_9x + s_{10}$\\
isomorphism:
&$\overline{s_1}\mapsto -\frac{3^6 5^2 7^2 11^2}{2^6 23^2 67^2}\overline{x}+\frac{3^6 5^27^3 11^2 \cdot 19 \cdot 163}{2^{12}23^3 67^3}\overline{z}$\\
&$\overline{s_2}\mapsto -\sqrt{-\frac{3^{13}5^{5}7^{5}11^{5}}{2^{15}23^{5}67^{5}}}\overline{y}$\\
&$\overline{s_3}\mapsto \frac{3^95^37^311^3}{2^{10} 23^3 67^3}\overline{z}+\frac{3^{11}5^4 7^3 11^3 13}{2^{10} 23^4 67^4}\overline{x}^2-\frac{3^95^37^311^341\cdot 307\cdot 587\cdot 32677569187}{2^{28}23^7 67^7}\overline{y}^2$\\

&$-\frac{3^9 5^4 7^3 11^3 71\cdot 1759 \cdot 516147191239}{2^{27} 23^7 67^7}\overline{x}\overline{z} -\frac{3^9 5^4 7^3 11^3 31 \cdot 2280560407042079}{2^{30} 23^7 67^7}\overline{z}^2$
\\ & \\
\hline
\\ & \\
name: &$Z_{11}$\\
equation: &$f=x^3y+xy^4+y^5$\\
deformation: & $F=f+s_{1}y^4 + s_2xy^3 + s_3y^3 + + s_4xy^2 + s_5y^2 + s_6xy + s_7y s_8x^2+ s_9x + s_{10}$\\
isomorphism:
&$\overline{s_1} \mapsto -\frac{2^{28} 3^3 7^4 11^4}{2399^4}\overline{x}-\frac{2^{29} 3^2 7^4 11^4 23 \cdot 53 \cdot 4405133 }{5^3 2399^6}\overline{y}^2$\\
&$\overline{s_2} \mapsto -\frac{2^{20} 3^2 7^3 11^3 173 \cdot 5879}{5^2 2399^4}\overline{x}+\frac{2^{20} 3^3 7^3 11^3}{2399^3}\overline{y}$\\
&$+\frac{2^{18} 7^3 11^3 59 \cdot 569 \cdot 49081 \cdot 52566671 \cdot 113887106221771273}{3^2 5^7 2399^7 271}\overline{x}\overline{y}$\\
&$+\frac{2^{18}7^311^341\cdot 13677187 \cdot 109919494930768288379}{3 \cdot 5^5 2399^6 271}\overline{y}^2$
\end{tabular}

 See appendix~\ref{sec:algorithm} for an algorithm that calculates these isomorphormism.
\end{exmp}

In the table in appendix~\ref{sec:appendix-table} we list the relevant equations of the modular strata of all 14 unimodal exceptional singularities.
The last monomial of the normal form induces in $Q(f)$ the Hesse form of the quasihomogeneous leading form $f_0$. 
The representatives of a $\C$-basis of $T(f)$ $b_1,\ldots, b_\tau$ are given 
in the second column, 
which are used to fix a miniversal deformation $F=f+\sum_i s_ib_i$. 
The equations in column three describe the modular stratum in the base of this deformation. 
Note, that we omit equations of the form $s_i +\CO(s^2)$, because we may eliminate 
the variables occurring linearly 
in these equations.
A minimal embedding with the remaining variables is given by the listed equations 
in the last column.
In the miniversal deformation these remaining variables $s_i$ 
correspond to the monomials printed in bold.

\section{Further examples and questions} 
We have calculated modular strata  for singularities of higher modality, too. 
The results raise hope that our observation generalises. 
We give one example of a bimodal singularity.
\begin{exmp}
\label{sec:exmp-complicated}
  We consider the hypersurface singularity given by the 
  semi-quasihomogeneous singularity of Hesse type
  $f=x^{10}+y^3+x^4y^2$. A miniversal deformation is defined  by 
  \begin{eqnarray*}
f&=&\mathbf{s_1}+\mathbf{s_2}x +\mathbf{s_3}x^2+ \mathbf{s_4}x^3 +\mathbf{s_5}x^4+ \mathbf{s_6}x^5 +\mathbf{s_7}x^6 +\mathbf{s_8}x^7 +\mathbf{s_9}x^8 +\mathbf{s_{10}}y\\
 &&+\mathbf{s_{11}}xy +\mathbf{s_{12}}x^2y+ \mathbf{s_{13}}x^3y+ \mathbf{s_{14}}x^4y +\mathbf{s_{15}}x^5y +\mathbf{s_{16}}x^6y +\mathbf{s_{17}}x^7y   
  \end{eqnarray*}
 The maximal modular subgerm $M$ in the base this deformation is defined by an ideal generated by
 \begin{eqnarray*}
   \mathbf{s_1} &+& \CO(\mathbf{s^2}),\\ &\vdots& \\ \mathbf{s_8} &+& \CO(\mathbf{s^2}), \\
   \mathbf{s_9^2}&-& \textstyle \frac{9}{256}\mathbf{s_{17}^4s_9}-\frac{29342801}{335104000}\mathbf{s_{17}^6s_9}
   -\frac{9963}{343146496}\mathbf{s_{17}^8}-\frac{831341932017399}{3283872972800000}\mathbf{s_{17}^{10}},\\
   \mathbf{s_{10}} &+& \CO(s^2),\\  &\vdots&\\ \mathbf{s_{16}} &+& \CO(\mathbf{s^2}),\\
   \mathbf{s_{17}^9} &-& \textstyle\frac{67372}{106029}\mathbf{s_{17}^7s_9}.
 \end{eqnarray*}
 The local ring $\CO_M = \CO_{17}/J_M$ is again isomorphic to $Q(f)$ via 
\begin{eqnarray*}
  \varphi:\CO_M &\rightarrow& \C\{x,y\}/(\textstyle \frac{\partial f}{\partial x}, \frac{\partial f}{\partial y}),\\
  \overline{\mathbf{s_{17}}} &\mapsto& a_1x,\\
  \overline{\mathbf{s_{9}}} &\mapsto& a_2by+a_3x^4+a_4x^2y+a_5x^6+a_6x^4y+a_7x^8,
\end{eqnarray*} 
with coefficients
\begin{eqnarray*}
  a_1&=& \textstyle 8\sqrt[4]{\frac{17943573032}{1269497754275}},\\
  a_2&=& \textstyle \frac{2261952}{84215}\sqrt{\frac{17943573032}{1269497754275}},\\
  a_3&=& \textstyle \frac{753984}{84215}\sqrt{\frac{17943573032}{1269497754275}}+\frac{1291937258304}{1269497754275},\\
  a_4&=& \textstyle \frac{9220238621242928785663198981632}{1007265342568292675484765625},\\
  a_5&=& \textstyle \frac{25742505984143872}{158687219284375}\sqrt{\frac{17943573032}{1269497754275}}+
  	\frac{3073412873747642928554399660544}{1007265342568292675484765625)},\\
  a_6&=&\textstyle \frac{547510092328050056695293440974819328}{377724503463109753306787109375}\sqrt{\frac{17943573032}{1269497754275}},\\
  a_7&=&\textstyle \frac{547510092328050056695293440974819328}{1133173510389329259920361328125}\sqrt{\frac{17943573032}{1269497754275}}.
\end{eqnarray*} 
\end{exmp}

In all the examples we have considered a function $f$ defining an isolated hypersurface singularity $X_0$, 
and related its modular stratum to the Milnor algebra of $f$. 
If we take another $\CK$-equivalent function $f'$, the isomorphism-class 
of the modular stratum does not change by definition. 
While $\mu(f)$ is an invariant of $\CK$-class, 
this is in general not true  for the isomorphism-class of the Milnor algebra, cf. \cite{by}.

Nevertheless, we have the following proposition  for singularities with $\tau=\mu-1$.

\begin{prop}
  Let be $f$ an analytic function with isolated critical point with $\tau(f)=\mu(f)-1$, 
  then its Milnor algebra is $\CK$-invariant.
\end{prop}
  This is a direct implication of the following two lemmas.  
  \begin{lem}\cite{by}
   \label{lemma:bensonyau}
    If $\mathfrak{m} f \subset \mathfrak{m}(\frac{\partial{f}}{\partial x_1},¸\ldots, 
    \frac{\partial{f}}{\partial x_n})$ then $\mathcal{RL}$-class and the $\CK$-class of $f$ coincide.
  \end{lem}
  
  \begin{lem}
    \label{lemma:saito}
    Let $f\in\C\{x_1,\ldots,x_n\}$ be a function with isolated singularity and $\delta \in \mathcal{T}_{\C^n}$ a vector field with
    $\delta(f) \in (f)$ then $\delta \in \mathfrak{m}\mathcal{T}_{\C^n}$.
  \end{lem}

  \begin{proof}
    This statement is only a slight generalisation of \cite[lemma~4.2]{saqh} 
    and we can use the same proof found there with
    $g f=\delta(f)$ instead of $f$.
  \end{proof}

\begin{proof}[Proof of the Proposition]  
  Look at the exact sequence 
  $$\begin{CD}
    0 \rightarrow \Ann(f) \rightarrow Q(f) @>{\cdot \overline{f}}>> Q(f) \rightarrow T(f) \rightarrow 0.
  \end{CD}$$
  Then $\Ann(f)$ has $\C$-dimension $\mu-1$ and  thus equals the maximal ideal $\m_{Q(f)}$. That means
  $\mathfrak{m} f \subset (\frac{\partial{f}}{\partial x_1},¸\ldots, \frac{\partial{f}}{\partial x_n})$. 
  Because of lemma~\ref{lemma:saito}, even
  $\mathfrak{m} f \subset \mathfrak{m}(\frac{\partial{f}}{\partial x_1},¸\ldots, \frac{\partial{f}}{\partial x_n})$ 
  holds. By lemma~\ref{lemma:bensonyau} we know
  that the $\CK$-class of $f$ equals its $\mathcal{RL}$-class, but $Q(f)$ is $\mathcal{RL}$-invariant.
\end{proof}

Due to this result  we can speak of the Milnor algebra 
of a hypersurface singularity in the case $\tau=\mu-1$. 
Hence, we can state the following conjecture, motivated by our examples.

\begin{hypo}
  Consider a hypersurface singularity $f$ with $\tau=\mu-1$. 
  Then the local ring of the modular stratum $\CO_{M(f)}$ is of Milnor type, 
  i.e. there exists a germ of an analytic function $f'$ such that $Q(f') \cong \CO_M$. 
  If $f$ has an Artinian modular stratum, then the local ring of the modular stratum 
  is isomorphic to the Milnor algebra of $f$ itself.
\end{hypo}

\begin{rem}
 We found the modular strata to be of Milnor type in all computed examples. 
So one could ask more generally: 
{\em For which singularities is the modular stratum of Milnor type?}
\end{rem}

Up to now we were not able to proof the hypothesis. There are some indications that it is strongly connected with the
restriction $\mu-\tau=1$.
\newpage 

\begin{appendix}
\section{Modular strata of exceptional unimodal singularities}
\label{sec:appendix-table}
\begin{tabular}[htb]{|p{3.67cm}|p{2.4cm}|p{6.6cm}|}\hline
 Singularity  & $T(f)$-basis & Equations of Modular Stratum \\ \hline 
$E_{12}$: $x^3+y^7+z^2 + xy^5$ &  
$\mathbf{y^6, xy^3},$
$y^{5}, y^{4}, y^{3},$
$xy^{2}, y^{2}, xy,$
$y, x, 1$
& 
$s_{2}^{2}+\frac{60}{49}s_{1}^{2}s_{2}+\frac{900}{2401}s_{1}^{4}+\frac{96235}{5176556}s_{1}^{3}s_{2}+\frac{13397975}{1141430598}s_{1}^{5}+\frac{3278640522559975}{275929299213140034}s_{1}^{6}+\frac{2597795351536401419351125}{533624932064451711054208176}s_{1}^{7}$,

$s_{1}^{4}s_{2}-\frac{24110}{5537}s_{1}^{6}-\frac{1593048109045}{1338513732171}s_{1}^{7}$ \\ \hline
$E_{13}$: $x^3+y^8+z^2 + y^8$ & 
$\mathbf{y^{7},xy^{3}},y^{6},$
$y^{5},y^{4},y^{3},xy^{2},$
$y^{2},xy,y,x,1$
& $s_{2}^{2}-\frac{173}{160}s_{1}^{3}s_{2}-\frac{3}{16}s_{1}^{5}+\frac{559677}{985600}s_{1}^{6}-\frac{1324866803129}{315707392000}s_{1}^{7}+\frac{78583901548081399293}{2528184795136000000}s_{1}^{8}$,
$s_{1}^{4}s_{2}+\frac{56981}{160160}s_{1}^{7}-\frac{22773678397}{10260490240}s_{1}^{8}$
\\ \hline
$E_{14}$: $x^3+xy^5+z^2 + xy^6$ & 
$\mathbf{y^7, xy^4},y^{6},y^{5},$
$y^{4},xy^{3},y^{3},xy^{2},$
$y^{2},xy,y,x,1$
&
$s_{2}^{2}+\frac{39}{64}s_{1}^{2}s_{4}+\frac{1521}{16384}s_{1}^{4}-\frac{12861}{4259840}s_{1}^{4}s_{2}-\frac{16047261}{15267266560}s_{1}^{6}-\frac{17484782484339}{14450719403592908800}s_{1}^{8}$
$s_{1}^{5}s_{2}+\frac{85371}{260032}s_{1}^{7}$
\\ \hline
$Z_{11}$: $x^3y+y^5+z^2 + xy^4$ & 
$\mathbf{y^4,xy^3},y^{3},xy^{2}$
$y^{2},xy,y,x^{2},x,1$
& 
	$s_1^2s_2+\frac{727093}{124740}s_{1}s_2^3-\frac{20939}{10395}s_2^4-
	\frac{5728952745569}{5480985787200}s_2^5$,

	$s_{1}^3+\frac{3419}{405}s_{1}s_2^3-\frac{80}{27}s_2^4
	-\frac{320173073863}{213544900800}s_2^5$
\\ \hline
$Z_{12}$:  $x^3y+xy^4+z^2 + x^2y^3$ & 
$\mathbf{y^{5},xy^3},y^{4},y^{3},$
$xy^{2},y^{2},xy,$
$y,x^{2},x,1$
& 

$s_{1}s_{2}^{2}+\frac{9}{4}s_{1}^{4}-\frac{67}{80}s_{1}^{3}s_{2}+\frac{367587}{22400}s_{1}^{5}+\frac{23086009519}{250880000}s_{1}^{6}$,

$s_{2}^{3}+27s_{1}^{3}s_{2}-\frac{36837}{560}s_{1}^{5}-\frac{2311958449}{6272000}s_{1}^{6}$
 \\ \hline

$Z_{13}$:  $x^3y+y^6+z^2 + xy^5$ & 
$\mathbf{y^{5},xy^{4}},y^{4},xy^{3},$
$y^{3},xy^{2},y^{2},xy,$
$y,x^{2},x,1$
& $s_{1}^{3}-\frac{745}{1512}s_{1}s_{2}^{4}+\frac{1}{18}s_{2}^{5}$,

$s_{1}^{2}s_{2}-\frac{92291}{589680}s_{1}s_{2}^{4}+\frac{619}{49140}s_{2}^{5}$
\\ \hline
$W_{12}$:  $x^4+y^5+z^2 + x^2y^3$ & 
$\mathbf{y^{4},xy^{3}},y^{3},xy^{2},$
$y^{2},x^{2}y,xy,$
$y,x^{2},x,1$
& 
	$s_{1}^4-\frac{30445}{7392}s_{1}^2s_{2}^2+\frac{4240139}{1897280}s_{1}^3s_{2}^2$,

	$s_{2}^3-\frac{2696}{48125}s_{1}^3s_{2}$
   \\ \hline
$W_{13}$:  $x^4+xy^4+z^2 + y^6$ & 
$\mathbf{y^5, x^2y^2}, y^{4},y^{3},$
$xy^{2},y^{2},x^{2}y,$
$xy,y,x^{2},x,1$
& $s_{2}^{3}+s_{1}^{4}+\frac{767}{90}s_{1}^{2}s_{2}^{2}-\frac{730357}{132300}s_{1}^{6}$,

$s_{1}^{3}s_{2}-\frac{9623}{8190}s_{1}^{5}$
\\ \hline   
$Q_{10}$:  $x^3+y^4+yz^2 + xy^3 $ & 
$\mathbf{xz, xy^2, y^3},z,$
$y^{2},xy,y,x,1$
& $s_{1}^{2}+\frac{254}{735}s_{2}^2s_{3}-\frac{8}{49}s_{2}^{3}+\frac{14007328}{1897810425}s_{2}^{4}$,

$s_{1}s_{3}-\frac{70}{47}s_{1}s_{2}$,

$s_{3}^{2}+\frac{271}{11340}s_{2}^2s_{3}+\frac{31}{945}s_{2}^{3}-\frac{43015222}{7320125925}s_{2}^{4}$
\\ \hline
$Q_{11}$:  $x^3+y^2z+xz^3 + z^5 $ & $\mathbf{z^4,xy,xz^2},z^{3},$
$z^{2},xz,z,y,x,1$
& $s_{3}^{2}-\frac{3}{4}s_{1}^{3}+\frac{229}{40}s_{1}^{2}s_{3}-\frac{28137}{5600}s_{1}^{4}+\frac{21341173563}{172480000}s_{1}^{5}$,

$s_{1}s_{2}-\frac{81}{55}s_{3}s_{2}$, 

$s_{2}^{2}+\frac{8}{25}s_{1}^{2}s_{3}-\frac{951}{875}s_{1}^{4}+\frac{775383309}{26950000}s_{1}^{5}$
 \\ \hline
$Q_{12}$:  $x^3+y^5+yz^2 + xy^4 $ & 
$\mathbf{xz,xy^3,y^4},z,$
$y^{3},xy^{2},y^{2},$
$xy,y,x,1$
& $s_{1}^{2}-\frac{65}{2688}s_{3}s_{2}^{3}+\frac{1}{320}s_{2}^{4}$,

$s_{1}s_{3}-\frac{12}{7}s_{1}s_{2}$,

$s_{3}^{2}+\frac{145}{1008}s_{3}s_{2}^{3}-\frac{29}{1260}s_{2}^{4}$\\ \hline
$S_{11}$:  $x^4+y^2z+xz^2 + x^3z$ & 
$\mathbf{x^2z,x^2y,x^3},xz$
$z,xy,y,x^{2},x,1$
& 
	$s_3^2-s_1^3+\frac{81}{28}s_1^2s_3-\frac{1893181}{1293600}s_1^3s_3$,

	$s_3s_2-\frac{52}{231}s_1^2s_2$,

	$s_2^2+\frac{4}{9}s_1s_3-\frac{1}{5}s_1^3+\frac{457}{756}s_1^2s_3-
	\frac{513482}{1819125}s_1^3s_3$	
	 \\ \hline
$S_{12}$:  $x^2y+y^2z+xz^3 + z^5 $ & $\mathbf{z^4, xz^2, yz^2},$
$z^{3},z^{2},yz,$
$xz,z,y,x,1$
& 
$s_{3}^{2}-3s_{1}s_{2}-\frac{31}{20}s_{1}^{2}s_{3}+\frac{81119}{32760}s_{1}^{4}$,

$s_{3}s_{2}-\frac{1}{6}s_{1}^{3}+\frac{4729}{1170}s_{1}^{2}s_{4}+\frac{217334189}{223587000}s_{1}^{5}$,

$s_{2}^{2}+\frac{2}{3}s_{1}^{2}s_{3}-\frac{15077}{24570}s_{1}^{4}$\\ \hline
$U_{12}$: $x^3+y^3+z^4 + xyz^2 $ & 
$\mathbf{z^3, yz^2, xz^2},$
$z^{2},yz,xz,z,$
$xy,y,x,1$
& $s_{2}^{2}+\frac{21}{320}s_{1}^{2}s_{3}$,

$s_{3}^{2}+\frac{21}{320}s_{1}^{2}s_{2}$,

$s_{1}^{3}-\frac{1024}{45}s_{1}s_{2}s_{3}$ \\ \hline
\end{tabular}

\section{Finding isomorphisms between Artinian local rings}
\label{sec:algorithm}
We have to 
look for isomorphisms between local rings of 
modular strata and algebras of Milnor type
 in order to test the hypothesis for the above examples. 
Using computer algebra and a direct approach
we check, whether two Artinian local algebras  $A:=\C\llbracket x_1, \ldots, x_m\rrbracket /I_A$ 
and $B:=\C\llbracket x_1,\ldots,x_n\rrbracket /I_B$ are 
isomorphic. 
A slightly generalised problem is considered. 
Let $A$ be arbitrary and $B$ Artinian. We will 
check if there is a surjective homomorphism $A \rightarrow B$. 
After eliminating variables we may assume that $n = \embdim B$.

Choose representatives of a $\C$-basis of 
$B$ $\{f_1, \ldots ,f_r\} \subset \C[x_1, \ldots x_n]$ and make an ansatz as follows: 
\begin{eqnarray*}
  \varphi:\C\llbracket x_1, \ldots x_m\rrbracket  &\rightarrow& \C\llbracket x_1, \ldots x_n\rrbracket \\
  x_i &\mapsto& \sum_j \alpha_{ij}f_j
\end{eqnarray*}
Our aim is to decide for which choice of $\alpha_{11}, \ldots \alpha_{nr} \in \C$ 
$\varphi$
defines a surjection such that $\varphi(I_A) \subset I_B$. 
Put $R:=\C[\alpha_{11},\ldots,\alpha_{nr}]$ and
\begin{eqnarray*}
  \tilde{\varphi}:\C\llbracket x_1, \ldots x_n\rrbracket  &\rightarrow& R\llbracket x_1, \ldots x_n\rrbracket \\
  x_i &\mapsto& \sum_j \alpha_{ij}f_j.
\end{eqnarray*}
Fix generators $a_1, \ldots,  a_k$ of $I_A$ and a local term ordering $>$ on $R\llbracket x\rrbracket $
and compute the reduced normal forms $\tilde{a}_i := \NF_> (\tilde{\varphi}(a_i),G_B)$ of $a_i$ 
with respect to a standard basis $G_B$ of $I_B R\llbracket x\rrbracket $. 
Consider the ideal $J \subset R$, generated by the coefficients of the $\tilde{a}_i$. 
Any $p \in \C^{n\cdot r}$ defines an evaluation homomorphism 
$$e_p:R \rightarrow \C,\ \alpha \mapsto p.$$
\begin{lem} With the notation of above:
  $e_p \circ \tilde{\varphi}$ defines surjection 
  $\varphi:\C\llbracket X\rrbracket  \rightarrow \C\llbracket X\rrbracket $ with  $\varphi(I_A) \subset I_B$  
  if and only if $p \in V(J) \setminus V(\minor_n(\frac{\partial 
  \tilde{\varphi}}{\partial x}|_{x=0}))$. 
\end{lem}
\begin{proof}
  Obviously $e_p \circ \tilde{\varphi}$ defines indeed a surjection
  $\C\llbracket X\rrbracket  \rightarrow \C\llbracket X\rrbracket $ iff not all $n$-minors of $(\frac{\partial 
  \tilde{\varphi}}{\partial x}|_{x=0})$ vanish.
  The key observation is: $e_p(\tilde{a}_i)$ is a reduced normal form for $e_p(a_i)$ 
  with respect to $G_B$. 
  Hence, all coefficients of $\tilde{a}_i$ have to be mapped 
  to zero by $e_p$ in order to have $(e_p\circ\tilde{\varphi})(a_i) \in I_B$. 
  Thus, all coefficients must belong to $\ker e_p = \m_p$. 
\end{proof}
For Artinian algebras $A$ and $B$ we first check that their $\C$-dimensions coincide.
If the algorithm leads to a homomorphism $\varphi:A \rightarrow B$, this will be a surjective homomorphism of vector spaces of same dimension. Thus, $\varphi$ is injective, too.
If the algorithm does not succeed in finding a surjective homomorphism, then $A$ and $B$ can not 
be isomorphic. 

\newpage

\begin{algorithm}
  \caption{findSurjection$(I_A,I_B)$}
\begin{algorithmic}
  \STATE \textbf{Input:} \ \\
    \begin{itemize}
    \item $I_A =(a_1,\ldots,a_k)$ ideal in $\C\llbracket x_1,\ldots,x_m\rrbracket $
    \item $I_B =(b_1,\ldots,b_l)$ ideal in $\C\llbracket x_1,\ldots,x_n\rrbracket $
    \end{itemize}
  
  \STATE \textbf{Output:} $\varphi: \C\llbracket X\rrbracket  \iso \C\llbracket X\rrbracket $ such that $\varphi(I_A) \subset I_B$ or $\varphi = \texttt{nil}$
 
  \STATE $G_B \gets $ standardBasis($b_1, \ldots, b_l$)
  \STATE $\{f_1,\ldots, f_r\} \gets $ $\C$-Basis($B$)
  \STATE $R := \C[\alpha_{11},\ldots, \alpha_{nr}]$
  \STATE $\tilde{\varphi}  \gets \big[\C\llbracket X\rrbracket  \rightarrow R\llbracket X\rrbracket , \quad  x_i \mapsto \sum \alpha_{ij} f_j \big]$
  \STATE $J \gets (0)$
  \FOR{$i=1$ to $k$} 
   \STATE $\tilde{a}_i \gets \NF(\tilde{\varphi}(a_i),G_B)$
   \STATE $J \gets J+($coefficients of $\tilde{a}_i)$
  \ENDFOR
  \IF{$\minor_n(\frac{\partial \tilde{\varphi}}{\partial x}|_{x=0}) \not\subset J$} 
   \STATE choose $p\in V(J)\setminus V(\minor_n(\frac{\partial \tilde{\varphi}}{\partial x}|_{x=0}))$
   \STATE $\varphi \gets \big[\C\llbracket X\rrbracket  \rightarrow \C\llbracket X\rrbracket , \quad  x_i \mapsto \sum p_{ij} f_j \big]$
  \ELSE 
  \STATE \textbf{return} \texttt{nil}
  \ENDIF
\end{algorithmic}
\end{algorithm}
\end{appendix}

\enlargethispage{3em} 

\end{document}